\documentclass[12pt]{article}   	
\usepackage{geometry}                		
\usepackage[parfill]{parskip}    		
\usepackage{graphicx}				
\usepackage[pdftex,pagebackref,letterpaper=true,colorlinks=true,pdfpagemode=none,urlcolor=blue,linkcolor=blue,citecolor=blue,pdfstartview=FitH]{hyperref}
\usepackage{amsmath,amsfonts}
\usepackage{color}
\usepackage{amsthm}
\usepackage{boolexpr}
\usepackage{amssymb,latexsym,pstricks,pst-plot}
\usepackage[english]{babel}
\usepackage{pgf}
\usepackage{tikz}
\usetikzlibrary{arrows,automata}
\usepackage[latin1]{inputenc}
\usepackage[UKenglish]{isodate}
\allowdisplaybreaks

\usepackage{amssymb}
\newtheorem{thm}{Theorem}

\title{Structure and asymptotics for Catalan numbers modulo primes using automata}
\author{Rob Burns}

\begin{document}
\maketitle
\begin{abstract}
Let $C_n$ be the $n$th Catalan number. We show that the asymptotic density of the set \mbox{$\{n: C_n \equiv 0 \mod p \}$} is $1$ for all primes $p$, We also show that if $n = p^k -1$ then $C_n \equiv -1 \mod p$. Finally we show that if $n \equiv \{ \frac{p+1}{2}, \frac{p+3}{2}, ..., p-2 \} \mod p$ then $p$ divides $C_n$. All results are obtained using the automata method of Rowland and Yassawi.
\end{abstract}

\section{Introduction}
The {\it Catalan numbers\/} are defined by
$$
C_n := \frac{1}{n+1}\binom{2n}{n}.
$$

There has been much work in recent years and also going back to Kummer \cite{Kum1852} on analysing the Catalan numbers modulo primes and prime powers.  Deutsch and Sagan \cite{Sagan2006} provided a complete characterisation of Catalan numbers modulo 3. A characterisation of the Catalan numbers modulo 2 dates back to Kummer. Eu, Liu and Yeh \cite{Eu2008} provided a complete characterisation of Catalan numbers modulo 4 and 8. This was extended by Liu and Yeh \cite{Liu2010} to a complete characterisation modulo 8, 16 and 64. This result was restated in a more compact form by Kauers, Krattenthaler and M\"uller in \cite{KKM2012} by representing the generating function of $C_n$ as a polynomial involving a special function. The polynomial for $C_n$ modulo 4096 was also calculated. A method for extracting the coefficients of the generating function (i.e. $C_n$ modulo a prime power) was provided. Given the complexity of the polynomials (the polynomial for the $4096$ case takes a page and a half to write down) the computation would need to be done by computer. Krattenthaler and M\"uller  \cite{KM2013} used a similar method to examine $C_n$ modulo powers of $3$. They wrote down the polynomial for the generating function of $C_n$ modulo $9$ and $27$ and thereby generalised the mod $3$ result of \cite{Sagan2006}. The article by Lin \cite{Lin:2010ab} discussed the possible values of the odd Catalan numbers modulo $2^k$ and Chen and Jiang \cite{Chen2013} dealt with the possible values of the Catalan numbers modulo prime powers. 
Rowland and Yassawi \cite{RY2013} investigated $C_n$ in the general setting of automatic sequences.  The values of $C_n$ (as well as other sequences) modulo prime powers can be computed via automata. Rowland and Yassawi provided algorithms for creating the relevant automata. They established a full characterisation of $C_n$ modulo $\{2, 4, 8, 16, 3, 5 \}$ in terms of automata. They also extended previous work by establishing forbidden residues for $C_n$ modulo $\{32, 64, 128, 256, 512 \}$. In theory the automata can be constructed for any prime power but computing power and memory quickly becomes a barrier.

We will use Rowland and Yassawi's automata to establish asymptotic densities of $C_n$ modulo primes. We will also make note of some structure results that appear from an examination of the relevant state diagrams of the automata. Asymptotic densities for $C_n$ modulo $2^k$ and $3$ are available from \cite{Burns:2016v1}. In this paper we will be discussing primes $p \geq 5$.

Here, the asymptotic density of a subset $S$ of $\mathbb{N}$ is defined to be
$$
\lim_{N \to \infty} \frac {1}{N} \#\{ n \in S : n \leq N \}
$$
if the limit exists, where $\#S$ is the number of elements in a set $S$. 

In particular, we will show that the asymptotic density of the set
\begin{equation}
S_p( \, 0 \,) = \{n: C_n \equiv 0 \mod p \}
\end{equation}
is $1$ for all primes $p$. For $p \in \{2 , 3 \}$ this result was shown in \cite{Burns:2016v1}. We will also show that if $n \equiv \{ \frac{p+1}{2}, \frac{p+3}{2}, ..., p-2 \} \mod p$ then $p$ divides $C_n$. Finally,
we will show that if \mbox{$n = p^k - 1$} for prime $p$ then \mbox{$C_n \equiv -1 \mod p$}. A stronger result has been known for $p=2$ since Kummer. Namely,

\bigskip
\begin{thm}
\label{mod2}
For all $n \geq 0$, $C_n$ is odd if and only if $n = 2^k - 1$ for some $k \geq 0$.
\end{thm}
\bigskip
The "only if" part does not apply for primes greater than $2$. For $p \in \{3, 5 \}$ the result is proved by Rowland and Yassawi in \cite{RY2013}. 
 
\bigskip

For a number $p$, we write the base $p$ expansion of a number $n$ as 
$$
[\, n ]\,_{p} = \langle n_{r} n_{r-1} ... n_{1} n_{0} \rangle
$$ 
where $n_{i} \in [\, 0, p - 1]\,$  and 
$$
n = n_{r} p^{r} + n_{r-1} p^{r-1} + ... + n_{1} p + n_{0}.
$$ 

\bigskip

\section{Background on automata for $C_n \mod p$}
Rowland and Yassawi showed in \cite{RY2013} that the behaviour of sequences such as $C_n~\mod~p$ can be studied by the use of finite state automata. The automaton has a finite number of states and rules for transitioning from one state to another. In the form described in \cite{RY2013} (algorithm 1) each state $s$ is represented by a polynomial in 2 variables $x$ and $y$. Each state has a value obtained by evaluating the polynomial at $x = 0$ and $y = 0$. All calculations are made modulo $p$. For the Catalan case the initial state
$s_1$ is represented by the polynomial
\begin{equation}
\label{R}
R(\, x, y \, )\,  = \, y(\, 1 - 2xy - 2xy^2 \, ).
\end{equation}

New states are constructed by applying the Cartier operator $\Lambda_{d, d}$ to the polynomials
$$
s_i*Q(\, x, y \, )^{p-1}
$$
for $d \in \{0, 1, ... , p-1 \}$ where $\{ s_i \}$ are the already calculated states and the polynomial $Q$ is defined by
\begin{equation}
\label{Q}
Q(\, x, y \, )\,  = \, x(\, y + 1 \,)^2 - 1.
\end{equation}
The Cartier operator is a linear map on polynomials defined by
$$
\Lambda_{d_1, d_2} (\, \sum_{m, n \geq 0} a_{m, n} \, x^{m} \, y^{n} \, ) = \sum_{m, n \geq 0} a_{pm + d_1, pn + d_2} \, x^{m} y^{n}.
$$

Since the Cartier operator maintains or reduces the degree of the polynomial and there are only finitely many polynomials modulo $p$ of each degree, all states of the automaton are obtained within a known finite time. It will be seen later that the automaton has at most $p+3$ states. If
$$
\Lambda_{d, d} (s*Q^{p-1}) = t
$$
for states (i.e. polynomials) $s$ and $t$ then the transition from state $s$ to state $t$ under the input $d$ is part of the automaton.

To calculate $C_n \mod p$, $n$ is first represented in base $p$. The base $p$ digits of $n$ are fed into the automaton starting with the least significant digit. The automaton starts at the initial state $s_1$ and transitions to a new state as each digit is fed into it. The value of the final state after all $n$'s digits have been used is equal to $C_n~mod~p$. Refer to \cite{RY2013} for more details. 

In the remainder of this article we will provide details of the automata for a general prime $p \geq 5$. We will provide the polynomials and values for the states and the transitions between states. States are listed as $s_1$, $s_2$, ... . Transitions, when provided, will be in the form $(\, s, j)\,  \to t$ which means that if the automaton is in state $s$ and receives digit $j$ then it will move to state $t$.  We will call a state $s$ a {\em \bf loop} state if all transitions from $s$ go to $s$ itself, i.e. \mbox{$(\, s, j)\,  \to s$} for all choices of $j$.

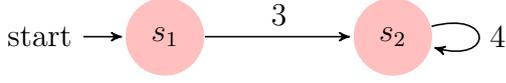
\begin{figure}[tbp]
\begin{tikzpicture}[->,>=stealth',shorten >=1pt,auto,node distance=3.0cm,
                    semithick]
  \tikzstyle{every state}=[fill=pink,draw=none,text=black]

  \node[initial,state] (A)                    {$s_1$};
  \node[state]         (B) [right of=A] {$s_2$};

  \path (A) edge    node {3} (B)
        (B) edge [loop right] node {4} (B);                      
 \end{tikzpicture}
                       \caption{Example of transition from state $s_1$ to state $s_2$ and a loop.}
                        \label{transition}
\end{figure}

States and transitions are represented visually in the form of a directed graph. For example, figure \ref{transition} represents an automaton which moves from state $s_1$ to state $s_2$ when it receives the digit $3$. It also moves from state $s_2$ to state $s_2$ (i.e. loops) if it is in state $s_2$ and receives a digit $4$.

\bigskip

\section{Preliminary calculations}

Before we start constructing the automata it will be convenient to first precompute $\Lambda_{d, d} (\, s(x, y)*Q(\, x, y \, )^{p-1} \,)$ for various choices of the polynomial $s$. The relevant results are contained in table \ref{lambdaactions}. When reading the table note that \mbox{$\binom{n}{m} = 0$} for \mbox{$m < 0$}. We will go through a few of the calculations from table~\ref{lambdaactions}.

\begin{table}[tbp]
  \centering
   \begin{tabular}{ | c ||  p{10cm} |}
 \hline
 & \\
State $s$ &  $\Lambda_{d, d} (\, s*Q^{p-1} \,)$ \\
& \\ \hline \hline
& \\
$1$ & $\binom{2d}{d}$ for $0 \leq d \leq p - 1$ \\
& \\
& \\ \hline
& \\
$y$ & $\binom{2d}{d-1}$ for $0 \leq d \leq p - 2$ \\
& \\
& $y + 1$ for $d = p - 1$ \\
& \\ \hline
& \\
$xy$ & $\binom{2d - 2}{d - 1}$ for $0 \leq d \leq p - 1$ \\
& \\
& \\ \hline
& \\
& $xy(\, y + 1 \,)$ for $d = 0$ \\
& \\
$xy^2$ & $\binom{2d - 2}{d - 2}$ for $1 \leq d \leq p - 1$ \\
& \\ \hline
& \\
$xy^3$ &  $-2xy(\, y + 1 \, )$ for $d = 0$ \\
& \\
&  $\binom{2d - 2}{d-3}$ for $1 \leq d \leq p - 2$ \\
& \\
& $y + 1$ for $d = p - 1$ \\
& \\ \hline
  \end{tabular}
  \caption{Table of values of $\Lambda_{d, d} (\, s*Q^{p-1} \,)$}
  \label{lambdaactions}
\end{table}

\bigskip

Firstly, the polynomial $Q^{p-1}$ can be written as
$$
Q^{p-1} (x, y) = (\, x(\, y + 1 \,)^2 - 1 \,)^{p-1}
$$
$$
= \sum_{k=0}^{p-1} \binom{p-1}{k} x^k (\, y+1 \,)^{2k} (\, -1 \,)^{p-1-k}
$$
\begin{equation}
\label{Qp-1}
= \sum_{k=0}^{p-1} \sum_{l=0}^{2k} \binom{p-1}{k} \binom{2k}{l} (\, -1 \,)^k x^k y^{l}.
\end{equation}
$$
= \sum_{k=0}^{p-1} \sum_{l=0}^{2k} a_{k, l} x^k y^l
$$
where
\begin{equation}
\label{akl}
a_{k, l} = \binom{p-1}{k} \binom{2k}{l} (\, -1 \,)^k
\end{equation}
for $0 \leq k \leq p-1$ and $0 \leq l \leq 2k$.

\bigskip

Then, for $r \geq 0$ and $t \geq 0$,

$$
\Lambda_{d, d} (\,x^r y^t Q(\, x, y \, )^{p-1} \,) = \Lambda_{d, d}  (\, \sum_{k=0}^{p-1} \sum_{l=0}^{2k} a_{k, l} x^{k+r} y^{l+t} \,)
$$
$$
=  \Lambda_{d, d}  (\, \sum_{k} \sum_{l} b_{k, l} x^{k} y^{l} \,)
$$
where 
\begin{equation}
\label{blk}
b_{k, l} = a_{k - r, l - t} = \binom{p-1}{k - r} \binom{2(k - r)}{l - t} (\, -1 \,)^{k - r}
\end{equation}
\bigskip
and the indices $k$ and $l$ in $b_{k, l}$ satisfy
\begin{equation}
\label{indexlimits}
r \leq k \leq p - 1 + r \, \; ; \, \;  t \leq l \leq 2(k - r) + t \leq 2p - 1 + t.
\end{equation}
\bigskip
So,
\begin{equation}
\label{lambdab}
\Lambda_{d, d} (\,x^r y^t Q(\, x, y \, )^{p-1} \,) = \sum_{k} \, \, \sum_{l} b_{pk + d, pl + d} x^k y^l.
\end{equation}

\bigskip
and the indices $k$ and $l$ in (\ref{lambdab}) satisfy
\begin{equation}
\label{indexlimitslambda}
r \leq pk + d \leq p - 1 + r \, \; ; \, \;  t \leq pl + d  \leq 2(pk + d - r) + t \leq 2p - 1 + t.
\end{equation}
\bigskip

We first compute $\Lambda_{d, d} (\,Q(\, x, y \, )^{p-1} \,)$ (i.e. $r=0$ and $t=0$ in (\ref{lambdab})). When $r=0$ and $t=0$ the only choice for $k$ satisfying the bounds in (\ref{indexlimitslambda}) is $k=0$. This then leaves $l=0$  as the only choice for $l$ satisfying (\ref{indexlimitslambda}). Then,
$$
\Lambda_{d, d} (\,Q(\, x, y \, )^{p-1} \,) = b_{d, d} = a_{d, d} = \binom{p-1}{d} \binom{2d}{d} (\, -1 \,)^d
$$
$$
\equiv \binom{2d}{d} \mod p
$$
since
$$
\binom{p-1}{d} \equiv (\, -1 \,)^d \mod p.
$$
\bigskip

This gives the first line of table~\ref{lambdaactions}. Note that \mbox{$\binom{0}{0} = 1$} and \mbox{$\binom{2d}{d} \equiv 0 \mod p$} for \mbox{$d \geq \frac{p+1}{2}$}.

\bigskip

We next calculate $\Lambda_{d, d} (\,y Q(\, x, y \, )^{p-1} \,)$, (i.e. $r=0$ and $t=1$ in (\ref{lambdab})). We have,
\bigskip
\begin{equation}
\label{bmn}
\Lambda_{d, d} (\,y Q(\, x, y \, )^{p-1} \,) = \sum_{k} \, \, \sum_{l} b_{pk+d, \, pl + d} \, x^k y^l.
\end{equation}
\bigskip
where $b_{k, l} = a_{k, l-1}$ for $0 \leq pk + d \leq p-1$ and \mbox{$1 \leq pl + d \leq 2(pk + d)+1$}.
If $d=0$, $k$ must be $0$ but there are then no choices for $l$ which fit the required bounds. Therefore,
$$
\Lambda_{0, 0} (\,yQ(\, x, y \, )^{p-1} \,)  = 0.
$$
\bigskip
If $1 \leq d \leq p-2$ then the only suitable choice for $\{k, l \}$ is $k=0$ and $l=0$. In this case
$$
\Lambda_{d, d} (\,y Q(\, x, y \, )^{p-1} \,)  = b_{d, d} = a_{d, d-1}
$$
\bigskip
$$
= \binom{p-1}{d} \binom{2d}{d - 1} (\, -1 \,)^d \equiv \binom{2d}{d-1} \mod p.
$$
\bigskip
We note that \mbox{$\binom{2d}{d-1} \equiv 0 \mod p$ for $d \geq \frac{p+1}{2}$}.
\bigskip

If $d=p-1$, there are two choices for $k$ and $l$ that satisfy the bounds in equation~(\ref{bmn}). These are $k=0$, $l=0$ and $k=0$, $l=1$. Therefore,
$$
\Lambda_{p-1, p-1} (\,y Q(\, x, y \, )^{p-1} \,)  = b_{p-1, p-1} + b_{p-1, 2p-1}y
$$
$$
= \binom{p-1}{p-1} \binom{2p-2}{p-2} (\, -1 \,)^{p-1} + \binom{p-1}{p-1} \binom{2p-2}{2p-2} (\, -1 \,)^{p-1} y.
$$
\bigskip
But $\binom{2p-2}{p-2} \equiv 1 \mod p$ so
$$
\Lambda_{p-1, p-1} (\,y Q(\, x, y \, )^{p-1} \,)  = y + 1.
$$
\bigskip

We next calculate $\Lambda_{d, d} (\,x \, y Q(\, x, y \, )^{p-1} \,)$, (i.e. $r=1$ and $t=1$ in (\ref{lambdab}) ). We have,
$$
\Lambda_{d, d} (\,x \, y Q(\, x, y \, )^{p-1} \,) = \sum_{k} \, \, \sum_{l} b_{pk+d, \, pl + d} \, x^k y^l.
$$
and from (\ref{indexlimitslambda}) $k$ and $l$ satisfy
$$
1 \leq pk + d \leq p \, \; ; \, \; 1 \leq pl + d \leq 2(pk + d - 1) + 1 \leq 2p.
$$
\bigskip
If $d=0$ there is only one choice for $\{k, l \}$ namely $k =1$ and $l=1$. In this case
$$
\Lambda_{0, 0} (\,x \, y Q(\, x, y \, )^{p-1} \,) = b_{p, p} \, x \, y
$$
$$
= \binom{p-1}{p-1} \binom{2p-2}{p-1} (\, -1 \,)^{p-1} \, x \, y  \equiv 0 \mod p.
$$
\bigskip
For $1 \leq d \leq p-1$ the only suitable choice for $\{ k, l \}$ is $k = 0$ and $l = 0$. So
$$
\Lambda_{d, d} (\,x \, y Q(\, x, y \, )^{p-1} \,) = b_{d, d}
$$
$$
= \binom{p-1}{d-1} \binom{2d-2}{d-1} (\, -1 \,)^{d-1} \equiv \binom{2d-2}{d-1} \mod p.
$$
\bigskip
We note that \mbox{$\binom{2d-2}{d-1} \equiv 0 \mod p$ for $d \geq \frac{p+3}{2}$}.
\bigskip

For brevity we will skip the calculation of $\Lambda_{d, d} (\,x \, y^2 Q(\, x, y \, )^{p-1} \,)$ and finally calculate $\Lambda_{d, d} (\,x \, y^3 Q(\, x, y \, )^{p-1} \,)$ (i.e. $r=1$ and $t=3$ in (\ref{lambdab})). We have,
$$
\Lambda_{d, d} (\,x \, y^3 Q(\, x, y \, )^{p-1} \,) = \sum_{k} \, \, \sum_{l} b_{pk+d, \, pl + d} \, x^k y^l.
$$
and from (\ref{indexlimitslambda}) $k$ and $l$ satisfy
\bigskip
\begin{equation}
\label{klxy3}
1 \leq pk + d \leq p \, \; ; \, \; 3 \leq pl + d \leq 2(pk + d - 1) + 3 \leq 2p.
\end{equation}
\bigskip
When $d = 0$, $k$ must be $1$ and $l$ can be $1$ or $2$. Therefore,

$$
\Lambda_{0, 0} (\,x \, y^3 Q(\, x, y \, )^{p-1} \,) = b_{p, p} \, x \, y  + b_{p, 2p} \, x \, y^2
$$
\bigskip
$$
= \binom{p-1}{p-1} \binom{2p-2}{p-3} (\, -1 \,)^{p-1} \, x \, y +  \binom{p-1}{p-1} \binom{2p-2}{2p-3} (\, -1 \,)^{p-1} \, x \, y^2
$$
\bigskip
$$
= -2xy - 2xy^2
$$
\bigskip 
since $\binom{2p-2}{p-3} \equiv -2 \mod p$ and $\binom{2p-2}{2p-3} \equiv -2 \mod p$.

\bigskip
For $d \in \{1, 2 \}$ there are no suitable choices for $k$ and $l$ in (\ref{klxy3}) so
$$
\Lambda_{d, d} (\,x \, y^3 Q(\, x, y \, )^{p-1} \,) = 0
$$
for $d \in \{1, 2 \}$.
For $3 \leq d \leq p-2$ the only suitable choice for $k$ in (\ref{klxy3}) is $k=0$. Then $l$ must also be $0$. Therefore, for $3 \leq d \leq p-2$,
$$
\Lambda_{d, d} (\,x \, y^3 Q(\, x, y \, )^{p-1} \,) = b_{d, d} = \binom{p-1}{d-1} \binom{2d-2}{d-3} (-1)^{d-1}
$$
$$
= \binom{2d-2}{d-3}
$$
since $\binom{p-1}{d-1} = (-1)^{d-1} \mod p$. Note that \mbox{$\binom{2d-2}{d-3} = 0$} for $d=1$ and $d=2$.

\bigskip
When $d = p-1$, $k=0$, $l=0$ and $k=0$, $l=1$ both satisfy the bounds in (\ref{klxy3}). So,
\bigskip
$$
\Lambda_{p-1, p-1} (\,x \, y^3 Q(\, x, y \, )^{p-1} \,) = b_{p-1, p-1} + b_{p-1, 2p-1}y
$$
\bigskip
$$
= \binom{p-1}{p-2} \binom{2p-4}{p-4} (-1)^{p-2} + \binom{p-1}{p-2} \binom{2p-4}{2p-4} (-1)^{p-2}y
$$
\bigskip
$$
\equiv y + 1 \mod p
$$
since $\binom{2p-4}{p-4} \equiv 1 \mod p$.

\bigskip

\section {Constructing the automata for $C_n \mod p$}

In this section we will describe the states and transitions of the automata for $C_n \mod p$. These are summarised in Table~\ref{results}. For given $d: 0 \leq d \leq p-1$ and state $s \in \{s_1, s_2, 1, -(y+1) \}$ the table gives the state equal to
$$
\Lambda_{d, d} (s*Q^{p-1}).
$$
The transition $(\, s, d \,) \to \Lambda_{d, d} (s*Q^{p-1})$ is then part of the automata.

\bigskip

\begin{table}[tbp]
  \centering
   \begin{tabular}{ | c || c | c | c | c |}
 \hline
 & & & &\\
 $d$ & $s_1$ & $s_2$ & $1$ & $-(y+1)$ \\ 
 & $y(\, 1 - 2xy - 2xy^2 \, )$ & $2xy(y+1)$ & & \\
 & & & &\\ \hline \hline
 & & & &\\
$d = 0$ & $s_2$& $s_2$ & $1$ & $-1$ \\ 
& & & &\\ \hline
& & & &\\
$1 \leq d \leq \frac{p-1}{2}$ & $\frac{1}{d} \binom{2d}{d-1}$ & $2 \binom{2d-1}{d}$ & $\binom{2d}{d}$ & $- \binom{2d+1}{d}$ \\ 
& & & &\\ \hline
& & & &\\
$\frac{p+1}{2} \leq d \leq p-2$ & $0$ & $0$ & $0$ & $0$ \\ 
& & & &\\ \hline
& & & &\\
$d=p-1$ & $-(y+1)$ & $0$ & $0$ & $-(y+1)$ \\
& & & &\\ \hline
  \end{tabular}
  \caption{Table of states and transitions.}
  \label{results}
\end{table}

\bigskip

\begin{figure}[tbp]
\begin{tikzpicture}[->,>=stealth',shorten >=1pt,auto,node distance=5.5cm,
                    semithick]
  \tikzstyle{every state}=[fill=pink,draw=none,text=black]

  \node[initial,state] (A)                    {$s_1$};
  \node[state]         (C) [below right of=A] {$C$};
  \node[state]         (B) [below of=C] {$-(y+1)$};
  \node[state]         (E) [right of=C]       {$0$};

  \path (A) edge [bend right] node {$p-1$} (B)
            edge node {$\{ \, 1, ... , \frac{p-1}{2} \, \}$} (C)
            edge [bend left] node {$\{ \, \frac{p+1}{2}, ... , p-2 \, \}$} (E)
        (B) edge [loop below] node {$p-1$} (B)
        edge node {$\{ \, 0, ... , \frac{p-3}{2} \, \}$} (C)
        edge [bend right] node  {$\{ \, \frac{p-1}{2}, ... , p-2 \, \}$} (E)
         (C) edge node  {$\{ \, \frac{p+1}{2}, ... , p-1 \, \}$} (E)
         edge [loop left] node  {$\{ \, 0, ... ,\frac{p-1}{2} \, \}$} (C)
                      (E) edge [loop right] node {all} (E); 
                      
 \end{tikzpicture}
                       \caption{Partial state diagram for $C_n \mod p$.}
                        \label{modpfigure1}
\end{figure}
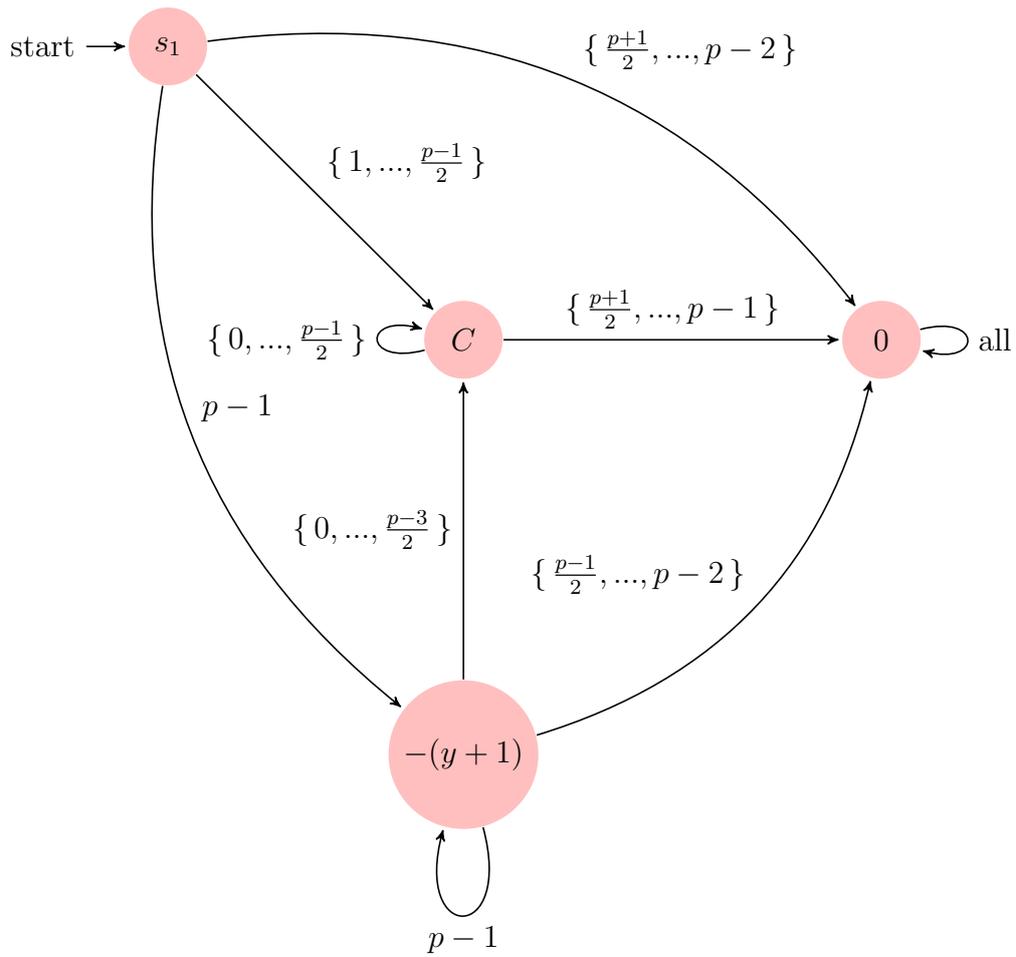

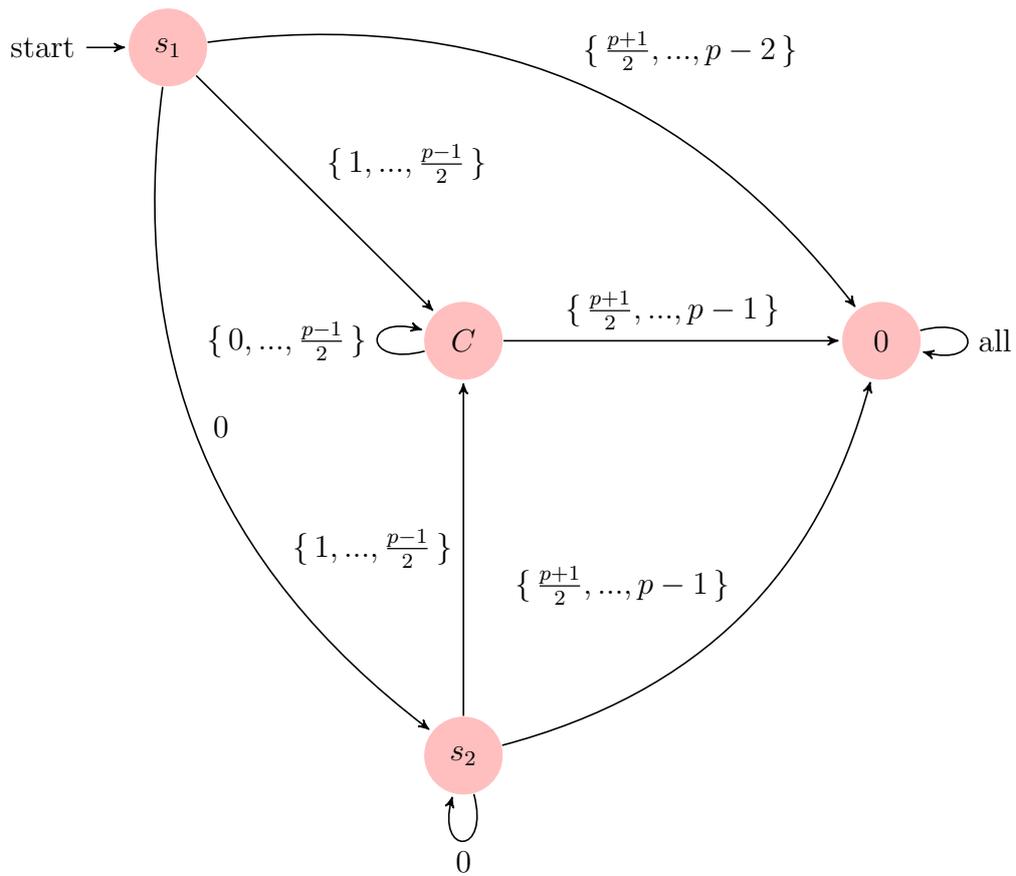
\begin{figure}[tbp]
\begin{tikzpicture}[->,>=stealth',shorten >=1pt,auto,node distance=5.5cm,
                    semithick]
  \tikzstyle{every state}=[fill=pink,draw=none,text=black]

  \node[initial,state] (A)                    {$s_1$};
  \node[state]         (C) [below right of=A] {$C$};
  \node[state]         (B) [below of=C] {$s_2$};
  \node[state]         (E) [right of=C]       {$0$};

  \path (A) edge [bend right] node {$0$} (B)
            edge node {$\{ \, 1, ... , \frac{p-1}{2} \, \}$} (C)
            edge [bend left] node {$\{ \, \frac{p+1}{2}, ... , p-2 \, \}$} (E)
        (B) edge [loop below] node {$0$} (D)
                      edge node  {$\{ \, 1, ... ,\frac{p-1}{2} \, \}$} (C)
                      edge [bend right] node  {$\{ \, \frac{p+1}{2}, ... , p-1 \, \}$} (E)
         (C) edge node  {$\{ \, \frac{p+1}{2}, ... , p-1 \, \}$} (E)
         edge [loop left] node  {$\{ \, 0, ... ,\frac{p-1}{2} \, \}$} (C)
                      (E) edge [loop right] node {all} (E); 
                      
 \end{tikzpicture}
                       \caption{Partial state diagram for $C_n \mod p$.}
                        \label{modpfigure2}
\end{figure}

Figures~\ref{modpfigure1} and \ref{modpfigure2} provide an alternative pictorial summary of the automata for $C_n~\mod~p$. We have broken the state diagram into 2 figures to improve clarity. In the figures the group of states labelled $C$ represents all states which consist of a constant non-zero polynomial. 

\bigskip

\bigskip

\bigskip

The calculation of the states will rely on the data contained in table~\ref{lambdaactions}. As mentioned earlier, the initial state $s_1$ for the automata is the polynomial defined in equation~\ref{R}. The second state $s_2$ is then given by
$$
s_2 = \Lambda_{0, 0} (\, s_1*Q(\, x, y \, )^{p-1} \,)
$$
$$
= \Lambda_{0, 0} (\, y(\, 1 - 2xy - 2xy^2 \,)*Q(\, x, y \, )^{p-1} \,)
$$
$$
= 0 -2(\, xy(\, y + 1 \,) ) -2(\, -2xy(\, y + 1 \,) \,)
$$
$$
= 2xy(\, y + 1 \,)
$$
\bigskip

For $1 \leq d \leq p -2$ we have
$$
\Lambda_{d, d} (\, s_1*Q(\, x, y \, )^{p-1} \,)
$$
$$
= \Lambda_{d, d} (\, y(\, 1 - 2xy - 2xy^2 \,)*Q(\, x, y \, )^{p-1} \,)
$$
$$
= \binom{2d}{d-1} -2 \binom{2d-2}{d-2} - 2 \binom{2d-2}{d-3}
$$
$$
= \binom{2d}{d-1} -2 \binom{2d-1}{d-2}
$$
$$
= \frac{(2d-1)!}{(d-1)!(d+1)!} (\, 2d - 2(\, d - 2 \,))
$$
$$
= \frac{1}{d}\binom{2d}{d-1}
$$
where we have used the identity 
\begin{equation}
\label{binid}
\binom{n}{m} + \binom{n}{m-1} = \binom{n+1}{m}.
\end{equation}
So these states are constant polynomials. Note that \mbox{$\binom{2d}{d-1} \equiv 0 \mod p$} for \mbox{$d \geq \frac{p+1}{2}$}. For $d=p-1$ we have
$$
= \Lambda_{p-1, p-1} (\, y(\, 1 - 2xy - 2xy^2 \,)*Q(\, x, y \, )^{p-1} \,)
$$
$$
= y + 1 - 2 \binom{2p-4}{p-3} - 2 (\, y + 1 \,)
$$
$$
= -(y + 1)
$$
since $\binom{2p-4}{p-3} \equiv 0 \mod p$.
\bigskip

The transitions for the state $s_2 = 2xy(y+1)$ can be worked out similarly. Using table~\ref{lambdaactions}, we have for $d=0$,
$$
\Lambda_{0, 0} (2xy(y+1)*Q^{p-1}) = 2xy(y+1) +2 \binom{-2}{-1} = s_2.
$$
\bigskip
For $1 \leq d \leq p - 1$,
$$
\Lambda_{d, d} (2xy(y+1)*Q^{p-1}) = 2 \binom{2d-2}{d-2} + 2 \binom{2d-2}{d-1}
$$
$$
= 2 \binom{2d-1}{d-1} = 2 \binom{2d-1}{d}
$$
where we have again used (\ref{binid}).

\bigskip

The transitions for the constant state with value $1$ can be taken straight from table~\ref{lambdaactions}, noting that $\binom{2d}{d} \equiv 0 \mod p$ for $d \geq \frac{p+1}{2}$. 

\bigskip

We finally consider the transitions for the state represented by the polynomial \mbox{$-(\, y + 1 \,)$}. For \mbox{$0 \leq d \leq p-2$},
$$
\Lambda_{d, d} (-(\, y + 1 \,) = - \binom{2d}{d-1} - \binom{2d}{d} = - \binom{2d+1}{d}.
$$
\bigskip
For $d = p-1$, 
$$
\Lambda_{p-1, p-1} (-(\, y + 1 \,) = -(y+1) - \binom{2p-2}{p-1} = -(y+1)
$$
since $\binom{2p-2}{p-1} \equiv 0 \mod p$.
\bigskip

Transitions for the constant states which have value other than $1$ can be obtained from the transitions for the state $1$ using the linearity of the Cartier operator.

\bigskip

\section{Conclusions}

Table \ref{results} shows that, under the algorithm we have used, the automata for $C_n \mod p$ contains at most $p+3$ states. The 3 polynomials \mbox{$y(\, 1 - 2xy - 2xy^2 \, )$}, \mbox{$2xy(\, y + 1 \,)$} and \mbox{$-(\, y + 1 \,)$} are always states of the automata. The other states are residues mod $p$. It is possible that all residues modulo $p$ always appear as states. In order to show that all residues appear it is sufficient to show that the set
$$
\{ \, \binom{2d}{d} :  \, 0 \leq d \leq \frac{p-1}{2} \} 
$$
generates $( \frac{ \mathbb{Z}}{p \mathbb{Z}} )^{\times}$. Forbidden residues do exist for $C_n \mod p^k$ for some primes $p$ and $k > 1$ as discussed in \cite{RY2013}. 

\bigskip

Table \ref{results} and figures \ref{modpfigure1} and \ref{modpfigure2} provide a clear view of the the working of $C_n \mod p$. Firstly, it is obvious that if $n = p^k - 1$ and so has a base $p$ representation in the form
$$
[\, n ]\,_{p} = \langle p-1, p-1, ..., p-1, p-1 \rangle
$$
then $C_n \equiv -1 \mod p$.

\bigskip

Secondly, the zero state is a loop state. Once the state path of a number reaches the zero state it cannot transition to any other state.

\bigskip

Thirdly, if the base $p$ representation of $n$ contains $1$ or more digits from the set \mbox{$\{\frac{p+1}{2}, \frac{p+3}{2}, ..., p-2 \}$} then \mbox{$C_n \equiv 0 \mod p$} since
$$
\Lambda_{d, d} (\, s Q^{p-1} \,) = 0
$$
for all states $s$ when $d \in \{ \frac{p+1}{2}, \frac{p+3}{2}, ..., p-2 \}$. Since the set of base-$p$ numbers which have no digits from the set $\{ \frac{p+1}{2}, \frac{p+3}{2}, ..., p-2 \}$ has asymptotic density $0$, it follows that the set $S_p(\, 0 \,)$ has asymptotic density $1$ for all $p$.

\bigskip

Fourthly, as a consequence of the point above, if $n \equiv \{ \frac{p+1}{2}, \frac{p+3}{2}, ..., p-2 \} \mod p$ then $p$ divides $C_n$. For $p=5$ this can be seen from the state diagram in \cite{RY2013}.

\bigskip

Finally, none of the constants contained within the group of states C in figures \ref{modpfigure1} and \ref{modpfigure2} can be 0. Therefore the criteria that $n$ must satisfy in order that $p \mid C_n$ can be deduced from the figures.

\bigskip

\section{Acknowledgement}
We would like to thank Eric Rowland for introducing us to automata as a tool for examining Catalan numbers.

\bigskip

\bibliographystyle{plain}
\begin{small}
\bibliography{ref}
\end{small}

\end{document}